# Dual Darboux Frame of a Spacelike Ruled Surface and Darboux Approach to Mannheim Offsets of Spacelike Ruled Surfaces

**Mehmet Önder[1], H. Hüseyin Uğurlu[2]**

[1]Celal Bayar University, Faculty of Arts and Sciences and, Department of Mathematics, Muradiye, Manisa, Turkey. E-mail: mehmet.onder@cbu.edu.tr

[2]Gazi University, Faculty of Education, Department of Secondary Education Science and Mathematics Teaching, Mathematics Teaching Program, Ankara, Turkey.

**Abstract**

In this paper, we define dual geodesic trihedron(dual Darboux frame) of a spacelike ruled surface. Then, we study Mannheim offsets of spacelike ruled surfaces in dual Lorentzian space by considering the E. Study Mapping. We represent spacelike ruled surfaces by dual Lorentzian unit spherical curves and define Mannheim offsets of the spacelike ruled surfaces by means of dual Darboux frame. We obtain relationships between the invariants of Mannheim spacelike offset surfaces and offset angle, offset distance. Furthermore, we give conditions for these surface offsets to be developable.

## 1. Introduction

In differential geometry, a surface is said to be "ruled" if through every point of the surface there is a straight line that lies on the surface. Then, a ruled surface can always be described (at least locally) as the set of points swept by a moving straight line. For example, a cone is formed by keeping one point of a line fixed whilst moving another point along a circle. Because of they are one of the simplest objects in geometric modeling, these surfaces are one of the most fascinating topics of the surface theory and also used in many areas of science such as Computer Aided Geometric Design (CAGD), mathematical physics, moving geometry, kinematics for modeling the problems and model-based manufacturing of mechanical products. For example, the building materials such as wood are straight and they can be considered as straight lines. The result is that if engineers are planning to construct something with curvature, they can use a ruled surface since all the lines are straight[9].

An offset surface is offset a specified distance from the original along the parent surface's normal. Offsetting of curves and surfaces is one of the most important geometric operations in CAD/CAM due to its immediate applications in geometric modeling, NC machining, and robot navigation[5]. Especially, the offsets of the ruled surfaces have an important role in (CAGD)[15,16]. The well-known offset of the ruled surfaces is Bertrand offsets which were defined by Ravani and Ku by considering a generalization of the theory of Bertrand curve for trajectory ruled surfaces in line geometry[17]. Moreover, there exists a one-to-one correspondence between the lines of line space and the points of dual unit sphere. This famous correspondence is known as E. Study Mapping[2]. Hence, the geometry of ruled surfaces can be thought as the geometry of dual spherical curves. By considering this correspondence, Küçük and Gürsoy have studied the integral invariants of closed Bertrand trajectory ruled surfaces in dual space[7]. They have given some characterizations of Bertrand offsets of trajectory ruled surfaces in terms of integral invariants (such as the angle of pitch and the pitch) of closed trajectory ruled surfaces and obtained the relationship between the area of projections of spherical images of Bertrand offsets of trajectory ruled surfaces and their integral invariants.

Recently, a new definition of special curve couple was given by Wang and Liu[22]. They have called these curves as Mannheim partner curves. Then, Orbay and et. al have given a generalization of the theory of Mannheim partner curves to the ruled surface and they have called this surface offset as Mannheim offset and shown that every developable ruled surface have a Mannheim offset if and only if an equation should be satisfied between the geodesic curvature and the arc-length of spherical indicatrix of the reference surface[9]. The corresponding characterizations of Mannheim offsets of ruled surfaces in Minkowski 3-space have been given in [10,14]. Furthermore, in [11] we have studied Mannheim offsets of ruled surfaces in dual space with Blaschke approach. We have given the characterizations of Mannheim offsets of ruled surfaces in terms of integral invariants of closed trajectory ruled surfaces and obtained the relationship between the area of projections of spherical images of Mannheim offsets of trajectory ruled surfaces and their integral invariants and we show that if the Mannheim offset surface are developable, then the striction lines of the surfaces are Mannheim Partner curves. Moreover, we have studied the Mannheim offsets of ruled surfaces in dual Lorentzian space by considering the Blaschke frame[12]. Also, we have given the dual Darboux frame of the timelike ruled surfaces with timelike rulings and study the Mannheim offsets of these surfaces[13].

In this paper, we define the dual geodesic trihedron(dual Darboux frame) of a spacelike ruled surface and give the real and dual curvatures of this surface. Then, we examine the Mannheim offsets of the spacelike ruled surfaces in view of their dual Darboux frame. Using the dual representations of the spacelike ruled surfaces, we give some theorems and new results which characterize the developable Mannheim spacelike surface offsets and we give the a new relationship between the developable Mannheim offsets and their striction lines.

## 2. Preliminaries

3-dimensional Minkowski space $IR_1^3$ is the real vector space $IR^3$ provided with the standard flat metric given by

$$\langle \vec{a}, \vec{a} \rangle = -a_1 b_1 + a_2 b_2 + a_3 b_3,$$

where $\vec{a}=(a_1,a_2,a_3)$ and $\vec{b}=(b_1,b_2,b_3) \in IR^3$. An arbitrary vector $\vec{a}=(a_1,a_2,a_3)$ of $IR_1^3$ is said to be timelike if $\langle \vec{a},\vec{a} \rangle < 0$, spacelike if $\langle \vec{a},\vec{a} \rangle > 0$ or $\vec{a}=0$, and lightlike (null) if $\langle \vec{a},\vec{a} \rangle = 0$ and $\vec{a} \neq 0$. Similarly, an arbitrary curve $\alpha(s)$ in $IR_1^3$ can locally be spacelike, timelike or null (lightlike), if all of its velocity vectors $\alpha'(s)$ are spacelike, timelike or null (lightlike), respectively[8]. The norm of a vector $\vec{a}$ is defined by $\|\vec{a}\| = \sqrt{|\langle \vec{a},\vec{a} \rangle|}$. Now, let $\vec{a}=(a_1,a_2,a_3)$ and $\vec{b}=(b_1,b_2,b_3)$ be two vectors in $IR_1^3$, then the Lorentzian cross product of $\vec{a}$ and $\vec{b}$ is given by

$$\vec{a} \times \vec{b} = \begin{vmatrix} \vec{e}_1 & -\vec{e}_2 & -\vec{e}_3 \\ a_1 & a_2 & a_3 \\ b_1 & b_2 & b_3 \end{vmatrix} = (a_2 b_3 - a_3 b_2, \ a_1 b_3 - a_3 b_1, \ a_2 b_1 - a_1 b_2)$$

where

$$\delta_{ij} = \begin{cases} 1 & i=j, \\ 0 & i \neq j, \end{cases} \quad \vec{e}_i = (\delta_{i1},\delta_{i2},\delta_{i3}) \text{ and } \vec{e}_1 \times \vec{e}_2 = -\vec{e}_3, \ \vec{e}_2 \times \vec{e}_3 = \vec{e}_1, \ \vec{e}_3 \times \vec{e}_1 = -\vec{e}_2.$$

By using this definition it can be easily shown that $\langle \vec{a} \times \vec{b}, \vec{c} \rangle = -\det(\vec{a},\vec{b},\vec{c})$ [20].

The sets of the unit timelike and spacelike vectors are called Lorentzian unit sphere and Lorentzian unit sphere, respectively, and denoted by

$$H_0^2 = \{\vec{a}=(a_1,a_2,a_3)\in E_1^3 : \langle \vec{a},\vec{a}\rangle = -1\}$$

and

$$S_1^2 = \{\vec{a}=(a_1,a_2,a_3)\in E_1^3 : \langle \vec{a},\vec{a}\rangle = 1\}$$

respectively(See [18]).

A surface in the Minkowski 3-space $IR_1^3$ is called a timelike surface if the induced metric on the surface is a Lorentz metric and is called a spacelike surface if the induced metric on the surface is a positive definite Riemannian metric, i.e., the normal vector on the spacelike (timelike) surface is a timelike (spacelike) vector[1].

## 3. Dual Numbers and Dual Lorentzian Vectors

Let $D = IR \times IR = \{\overline{a}=(a,a^*): a,a^* \in IR\}$ be the set of the pairs $(a,a^*)$. For $\overline{a}=(a,a^*)$, $\overline{b}=(b,b^*)\in D$ the following operations are defined on $D$:

Equality: $\overline{a}=\overline{b} \Leftrightarrow a=b,\ a^*=b^*$

Addition: $\overline{a}+\overline{b}=(a+b,\ a^*+b^*)$

Multiplication: $\overline{ab}=(ab,\ ab^*+a^*b)$

The element $\varepsilon = (0,1)\in D$ satisfies the relationships

$$\varepsilon \neq 0,\ \varepsilon^2 = 0,\ \varepsilon 1 = 1\varepsilon = \varepsilon. \tag{1}$$

Let consider the element $\overline{a}\in D$ of the form $\overline{a}=(a,0)$. Then the mapping $f:D\to IR$, $f(a,0)=a$ is a isomorphism. So, we can write $a=(a,0)$. Then, by the multiplication rule we have that

$$\overline{a} = (a,a^*)$$
$$= (a,0)+(0,a^*)$$
$$= (a,0)+(0,1)(a^*,0)$$
$$= a+\varepsilon a^*$$

Then $\overline{a}=a+\varepsilon a^*$ is called dual number and $\varepsilon$ is called dual unit. Thus the set of dual numbers is given by

$$D = \{\overline{a}=a+\varepsilon a^* : a,a^* \in IR,\ \varepsilon^2 = 0\}. \tag{2}$$

The set $D$ forms a commutative group under addition. The associative laws hold for multiplication. Dual numbers are distributive and form a ring over the real number field[2,4].

Dual function of dual number presents a mapping of a dual numbers space on itself. Properties of dual functions were thoroughly investigated by Dimentberg[3]. He derived the general expression for dual analytic (differentiable) function as follows

$$f(\overline{x}) = f(x+\varepsilon x^*) = f(x)+\varepsilon x^* f'(x), \tag{3}$$

where $f'(x)$ is derivative of $f(x)$ with respect to $x$ and $x, x^* \in IR$. This definition allows us to write the dual forms of some well-known functions as follows

$$\begin{cases} \cosh(\overline{x}) = \cosh(x+\varepsilon x^*) = \cosh(x)+\varepsilon x^* \sinh(x), \\ \sinh(\overline{x}) = \sinh(x+\varepsilon x^*) = \sinh(x)+\varepsilon x^* \cosh(x), \\ \sqrt{\overline{x}} = \sqrt{x+\varepsilon x^*} = \sqrt{x}+\varepsilon \dfrac{x^*}{2\sqrt{x}},\ (x>0). \end{cases} \tag{4}$$

Let $D^3 = D\times D\times D$ be the set of all triples of dual numbers, i.e.,

$$D^3 = \{\tilde{a} = (\bar{a}_1, \bar{a}_2, \bar{a}_3) : \bar{a}_i \in D, \ i = 1, 2, 3\}, \tag{5}$$

Then the set $D^3$ is called dual space. The elements of $D^3$ are called dual vectors. Analogue to the dual numbers, a dual vector $\tilde{a}$ may be expressed in the form $\tilde{a} = \vec{a} + \varepsilon \vec{a}^* = (\vec{a}, \vec{a}^*)$, where $\vec{a}$ and $\vec{a}^*$ are the vectors of $IR^3$. Then for any vectors $\tilde{a} = \vec{a} + \varepsilon \vec{a}^*$ and $\tilde{b} = \vec{b} + \varepsilon \vec{b}^*$ of $D^3$, the scalar product and the vector product are defined by

$$\langle \tilde{a}, \tilde{b} \rangle = \langle \vec{a}, \vec{b} \rangle + \varepsilon \left( \langle \vec{a}, \vec{b}^* \rangle + \langle \vec{a}^*, \vec{b} \rangle \right), \tag{6}$$

and

$$\tilde{a} \times \tilde{b} = \vec{a} \times \vec{b} + \varepsilon \left( \vec{a} \times \vec{b}^* + \vec{a}^* \times \vec{b} \right), \tag{7}$$

respectively, where $\langle \vec{a}, \vec{b} \rangle$ and $\vec{a} \times \vec{b}$ are the inner product and the vector product of the vectors $\vec{a}$ and $\vec{a}^*$ in $IR^3$, respectively.

The norm of a dual vector $\tilde{a}$ is given by

$$\|\tilde{a}\| = \|\vec{a}\| + \varepsilon \frac{\langle \vec{a}, \vec{a}^* \rangle}{\|\vec{a}\|}, \quad (\vec{a} \neq 0). \tag{8}$$

A dual vector $\tilde{a}$ with norm $1 + \varepsilon 0$ is called dual unit vector. The set of dual unit vectors is given by

$$\tilde{S}^2 = \{\tilde{a} = (a_1, a_2, a_3) \in D^3 : \langle \tilde{a}, \tilde{a} \rangle = 1 + \varepsilon 0\}, \tag{9}$$

and called dual unit sphere[2,4].

The Lorentzian inner product of two dual vectors $\tilde{a} = \vec{a} + \varepsilon \vec{a}^*$, $\tilde{b} = \vec{b} + \varepsilon \vec{b}^* \in D^3$ is defined by

$$\langle \tilde{a}, \tilde{b} \rangle = \langle \vec{a}, \vec{b} \rangle + \varepsilon \left( \langle \vec{a}, \vec{b}^* \rangle + \langle \vec{a}^*, \vec{b} \rangle \right),$$

where $\langle \vec{a}, \vec{b} \rangle$ is the Lorentzian inner product of the vectors $\vec{a}$ and $\vec{b}$ in the Minkowski 3-space $IR_1^3$. Then a dual vector $\tilde{a} = \vec{a} + \varepsilon \vec{a}^*$ is said to be timelike if $\vec{a}$ is timelike, spacelike if $\vec{a}$ is spacelike or $\vec{a} = 0$ and lightlike (null) if $\vec{a}$ is lightlike (null) and $\vec{a} \neq 0$ [18].

The set of all dual Lorentzian vectors is called dual Lorentzian space and it is denoted by

$$D_1^3 = \{\tilde{a} = \vec{a} + \varepsilon \vec{a}^* : \vec{a}, \vec{a}^* \in IR_1^3\}.$$

The Lorentzian cross product of dual vectors $\tilde{a}, \tilde{b} \in D_1^3$ is defined by

$$\tilde{a} \times \tilde{b} = \vec{a} \times \vec{b} + \varepsilon (\vec{a}^* \times \vec{b} + \vec{a} \times \vec{b}^*),$$

where $\vec{a} \times \vec{b}$ is the Lorentzian cross product in $IR_1^3$.

Let $\tilde{a} = \vec{a} + \varepsilon \vec{a}^* \in D_1^3$. Then $\tilde{a}$ is said to be dual unit timelike (resp. spacelike) vector if the vectors $\vec{a}$ and $\vec{a}^*$ satisfy the following equations:

$$<\vec{a}, \vec{a}> = -1 \ (resp. <\vec{a}, \vec{a}> = 1), \ <\vec{a}, \vec{a}^*> = 0. \tag{10}$$

The set of all unit dual timelike vectors is called the dual hyperbolic unit sphere, and is denoted by $\tilde{H}_0^2$,

$$\tilde{H}_0^2 = \{\tilde{a} = (a_1, a_2, a_3) \in D_1^3 : \langle \tilde{a}, \tilde{a} \rangle = -1 + \varepsilon 0\} \tag{11}$$

and the set of all unit dual spacelike vectors is called the dual Lorentzian unit sphere, and is denoted by $\tilde{S}_1^2$,

$$\tilde{S}_1^2 = \{\tilde{a} = (a_1, a_2, a_3) \in D_1^3 : \langle \tilde{a}, \tilde{a} \rangle = 1 + \varepsilon 0\}. \tag{12}$$

(For details see [18]).

**Definition 3.1.** *i) Dual Lorentzian timelike angle:* Let $\tilde{x}$ be a dual spacelike vector and $\tilde{y}$ be a dual timelike vector in $D_1^3$. Then the dual angle between $\tilde{x}$ and $\tilde{y}$ is defined by $<\tilde{x},\tilde{y}>=\|\tilde{x}\|\|\tilde{y}\|\sinh\bar{\theta}$. The dual number $\bar{\theta}=\theta+\varepsilon\theta^*$ is called the dual *Lorentzian timelike angle*.

*ii) Dual Central angle:* Let $\tilde{x}$ and $\tilde{y}$ be dual spacelike vectors in $D_1^3$ that span a dual timelike vector subspace. Then the dual angle between $\tilde{x}$ and $\tilde{y}$ is defined by $<\tilde{x},\tilde{y}>=\|\tilde{x}\|\|\tilde{y}\|\cosh\bar{\theta}$. The dual number $\bar{\theta}=\theta+\varepsilon\theta^*$ is called the dual *central angle*[18].

## 4. Dual Representation and Dual Darboux Frame of a Spacelike Ruled Surface

From E. Study mapping, the lines of the line space $IR^3$ correspond to the dual points(dual unit vectors) of the dual unit sphere[2,4]. Then, the dual spherical curve lying fully on $\tilde{S}^2$ represents a ruled surface of $IR^3$. In this section, we introduce this correspondence rule for spacelike ruled surfaces and give the dual geodesic trihedron(dual Darboux) frame of the spacelike ruled surfaces.

In the Minkowski 3-space $IR_1^3$, an oriented spacelike line $L$ is determined by a point $p\in L$ and a unit spacelike vector $\vec{a}$. Then, one can define $\vec{a}^*=\vec{p}\times\vec{a}$ which is called moment vector. The value of $\vec{a}^*$ does not depend on the point $p$, because any other point $q$ in $L$ can be given by $\vec{q}=\vec{p}+\lambda\vec{a}$ and then $\vec{a}^*=\vec{p}\times\vec{a}=\vec{q}\times\vec{a}$. Reciprocally, when such a pair $(\vec{a},\vec{a}^*)$ is given, one recovers the spacelike line $L$ as $L=\left\{(\vec{a}\times\vec{a}^*)+\lambda\vec{a}:\vec{a},\vec{a}^*\in E^3,\lambda\in IR\right\}$, written in parametric equations. The vectors $\vec{a}$ and $\vec{a}^*$ are not independent of one another and they satisfy the following relationships

$$\langle\vec{a},\vec{a}\rangle=1,\quad \langle\vec{a},\vec{a}^*\rangle=0 \tag{13}$$

The components $a_i$, $a_i^*$ ($1\leq i\leq 3$) of the vectors $\vec{a}$ and $\vec{a}^*$ are called the normalized Plucker coordinates of the spacelike line $L$. From (10), (11) and (13) we see that the dual spacelike unit vector $\tilde{a}=\vec{a}+\varepsilon\vec{a}^*$ corresponds to spacelike line $L$. This correspondence is known as E. Study Mapping: There exists a one-to-one correspondence between the spacelike vectors of dual Lorentzian unit sphere $\tilde{S}_1^2$ and the directed spacelike lines of the Minkowski space $IR_1^3$[18]. By the aid of this correspondence, the properties of the spatial motion of a spacelike line can be derived. Hence, the geometry of spacelike ruled surface is represented by the geometry of dual Lorentzian curve lying on the dual Lorentzian unit sphere $\tilde{S}_1^2$.

Veldkamp introduced the dual representation and dual geodesic trihedron of a ruled surface[21]. Now, we use the similar procedure to give the dual Darboux frame of a spacelike ruled surface.

Let $(\tilde{k})$ be a dual Lorentzian curve represented by the dual spacelike unit vector $\tilde{e}(u)=\vec{e}(u)+\varepsilon\vec{e}^*(u)$. The real unit vector $\vec{e}$ draws a curve on the real Lorentzian unit sphere $S_1^2$ and is called the (real) indicatrix of $(\tilde{k})$. We suppose throughout that it is not a single point. We take the parameter $u$ as the arc-length parameter $s$ of the real indicatrix and denote the differentiation with respect to $s$ by primes. Then we have $\langle\vec{e}',\vec{e}'\rangle=-1$. The vector $\vec{e}'=\vec{t}$ is the unit vector parallel to the tangent of the indicatrix and it is also the unit normal of the surface. The equation $\vec{e}^*(s)=\vec{p}(s)\times\vec{e}(s)$ has infinity of solutions for the function $\vec{p}(s)$. If we take $\vec{p}_o(s)$ as a solution, the set of all solutions is given by $\vec{p}(s) = \vec{p}_o(s) + \lambda(s)\vec{e}(s)$, where $\lambda$ is a real scalar function of $s$. Therefore we have $\langle \vec{p}', \vec{e}' \rangle = \langle \vec{p}'_o, \vec{e}' \rangle + \lambda$. By taking $\lambda = \lambda_o = -\langle \vec{p}'_o, \vec{e}' \rangle$ we see that $\vec{p}_o(s) + \lambda_o(s)\vec{e}(s) = \vec{c}(s)$ is the unique solution for $\vec{p}(s)$ with $\langle \vec{c}', \vec{e}' \rangle = 0$. Then, the given dual curve $(\tilde{k})$ corresponding to the spacelike ruled surface

$$\varphi_e = \vec{c}(s) + v\vec{e}(s) \tag{14}$$

may be represented by

$$\tilde{e}(s) = \vec{e} + \varepsilon \vec{c} \times \vec{e} \tag{15}$$

where

$$\langle \vec{e}, \vec{e} \rangle = 1, \quad \langle \vec{e}', \vec{e}' \rangle = -1, \quad \langle \vec{c}', \vec{e}' \rangle = 0.$$

and $\vec{c}$ is the direction vector of the striction curve. Then we have

$$\|\tilde{e}'\| = 1 + \varepsilon \det(\vec{c}', \vec{e}, \vec{t}) = 1 + \varepsilon \Delta \tag{16}$$

where $\Delta = \det(\vec{c}', \vec{e}, \vec{t})$ which characterizes the developable spacelike surface, i.e, the spacelike surface is developable if and only if $\Delta = 0$. Then, the dual arc-length $\bar{s}$ of the dual curve $(\tilde{k})$ is given by

$$\bar{s} = \int_0^s \|\tilde{e}'(u)\| du = \int_0^s (1 + \varepsilon \Delta) du = s + \varepsilon \int_0^s \Delta du \tag{17}$$

From (17) we have $\bar{s}' = 1 + \varepsilon \Delta$. Therefore, the dual unit tangent to the dual curve $\tilde{e}(s)$ is given by

$$\frac{d\tilde{e}}{d\bar{s}} = \frac{\tilde{e}'}{\bar{s}'} = \frac{\tilde{e}'}{1+\varepsilon\Delta} = \tilde{t} = \vec{t} + \varepsilon(\vec{c} \times \vec{t}) \tag{18}$$

Introducing the dual unit vector $\tilde{g} = -\tilde{e} \times \tilde{t} = \vec{g} + \varepsilon \vec{c} \times \vec{g}$ we have the dual frame $\{\tilde{e}, \tilde{t}, \tilde{g}\}$ which is known as dual geodesic trihedron or dual Darboux frame of $\varphi_e$ (or $(\tilde{e})$). Also, it is well known that the real orthonormal frame $\{\vec{e}, \vec{t}, \vec{g}\}$ is called the geodesic trihedron of the indicatrix $\vec{e}(s)$ with the derivations

$$\vec{e}' = \vec{t}, \quad \vec{t}' = \vec{e} + \gamma \vec{g}, \quad \vec{g}' = \gamma \vec{t} \tag{19}$$

where $\gamma$ is called the conical curvature[6,19].

Lot now consider the derivations of vectors of dual geodesic trihedron and find the dual Darboux formulae of a spacelike ruled surface.

From (18) we have $\langle \tilde{t}, \tilde{t} \rangle = -1 + \varepsilon 0$. By using this equality and considering that $\tilde{g} = -\tilde{e} \times \tilde{t}$, we have

$$\left\langle \tilde{t}, \frac{d\tilde{t}}{d\bar{s}} \right\rangle = 0, \quad \frac{d\tilde{g}}{d\bar{s}} = -\tilde{e} \times \frac{d\tilde{t}}{d\bar{s}}. \tag{20}$$

For the derivative of $\tilde{t}$ let write

$$\frac{d\tilde{t}}{d\bar{s}} = \bar{\alpha}\tilde{e} + \bar{\beta}\tilde{t} + \bar{\gamma}\tilde{g} \tag{21}$$

where $\bar{\alpha}, \bar{\beta}, \bar{\gamma}$ are the dual functions of dual arc-length $\bar{s}$. The first equation of (20) gives that $\bar{\beta} = 0$. Thus from the second equation of (20) we have

$$\frac{d\tilde{g}}{d\bar{s}} = \bar{\gamma}\tilde{t}. \tag{22}$$

Finally, from (22) and the equality $\tilde{t} = \tilde{g} \times \tilde{e}$, we obtain

$$\frac{d\tilde{t}}{ds} = \tilde{e} + \overline{\gamma}\tilde{g}. \tag{23}$$

Then from (18), (22) and (23) we have the following theorem.

**Theorem 4.1.** *The derivatives of the vectors of dual frame* $\{\tilde{e},\tilde{t},\tilde{g}\}$ *of a spacelike ruled surface are obtained as follows*

$$\frac{d\tilde{e}}{ds} = \tilde{t}, \quad \frac{d\tilde{t}}{ds} = \tilde{e} + \overline{\gamma}\tilde{g}, \quad \frac{d\tilde{g}}{ds} = \overline{\gamma}\tilde{t} \tag{24}$$

*and called dual Darboux formulae of the spacelike ruled surface* $\varphi_e$. *Then the dual Darboux vector of the trihedron is* $\tilde{d} = -\overline{\gamma}\tilde{e} + \tilde{g}$.

Let now give the invariants of the surface. Since $\tilde{s}' = 1 + \varepsilon\Delta$, it follows from (22) that
$$\tilde{g}' = \overline{\gamma}(1+\varepsilon\Delta)\tilde{t} \tag{25}$$
On the other hand using that $\tilde{g} = \vec{g} + \varepsilon\vec{c}\times\vec{g}$, from (18) we have
$$\begin{aligned}\tilde{g}' &= \gamma\vec{t} + \varepsilon(\vec{c}'\times\vec{g}\gamma\vec{c}\times\vec{t}) \\ &= \gamma\vec{t} + \varepsilon\gamma(\vec{c}\times\vec{t}) + \varepsilon(\vec{c}'\times\vec{g}) \\ &= \gamma\tilde{t} + \varepsilon(\vec{c}'\times\vec{g})\end{aligned} \tag{26}$$
Therefore, from (25) and (26) we obtain
$$\overline{\gamma}(1+\varepsilon\Delta)\tilde{t} = \gamma\tilde{t} + \varepsilon(\vec{c}'\times\vec{g}) \tag{27}$$
which gives us
$$\overline{\gamma}(1+\varepsilon\Delta) = \gamma - \varepsilon\delta \tag{28}$$
where $\delta = \langle \vec{c}',\vec{e}\rangle$. Then from (28) we have
$$\overline{\gamma} = \gamma - \varepsilon(\delta + \gamma\Delta) \tag{29}$$
Moreover, since $\vec{c}'$ as well as $\vec{e}$ is perpendicular to $\vec{t}$, for the real scalar $\mu$ we may write $\vec{c}'\times\vec{e} = \mu\vec{t}$. Then
$$\Delta = \det(\vec{c}',\vec{e},\vec{t}) = -\langle \vec{c}'\times\vec{e},\vec{t}\rangle = -\mu\langle \vec{t},\vec{t}\rangle = \mu.$$
Hence $\vec{e}\times(\vec{c}'\times\vec{e}) = \Delta\vec{e}\times\vec{t} = -\Delta\vec{g}$ and $\vec{c}' = \delta\vec{e} - \Delta\vec{g}$.

The functions $\gamma(s)$, $\delta(s)$ and $\Delta(s)$ are the invariants of the spacelike ruled surface $\varphi_e$. They determine the spacelike ruled surface uniquely up to its position in the space. For example, if $\delta = \Delta = 0$ we have that $\vec{c}$ is constant. It means that the spacelike ruled surface $\varphi_e$ is a spacelike cone.

### 4.1. Elements of Curvature of a Dual Lorentzian Curve

The dual radius of curvature of dual Lorentzian curve(spacelike ruled surface) $\tilde{e}(s)$ is can be calculated analogous to common Lorentzian differential geometry of curves as follows

$$\overline{R} = \frac{\left\|\frac{d\tilde{e}}{ds}\right\|^3}{\left\|\frac{d\tilde{e}}{ds}\times\frac{d^2\tilde{e}}{ds^2}\right\|} = \frac{1}{\sqrt{1+\overline{\gamma}^2}} \tag{30}$$

The unit vector $\tilde{d}_o$ with the same sense as the Darboux vector $\tilde{d} = -\overline{\gamma}\tilde{e} + \tilde{g}$ is given by

$$\tilde{d}_o = -\frac{\overline{\gamma}}{\sqrt{1+\overline{\gamma}^2}}\tilde{e} + \frac{1}{\sqrt{1+\overline{\gamma}^2}}\tilde{g} \tag{31}$$

It is clear that $\tilde{d}_o$ is spacelike. Then, the dual angle between $\tilde{d}_o$ and $\tilde{e}$ satisfies the followings

$$\sin \bar{\rho} = \frac{1}{\sqrt{1+\bar{\gamma}^2}} = \bar{R}, \quad \cos \bar{\rho} = \frac{-\bar{\gamma}}{\sqrt{1+\bar{\gamma}^2}}$$

where $\bar{\rho}$ is the dual spherical radius of curvature. Hence

Moreover, in [13], we have given the invariants, dual radius of curvature and dual spherical radius of curvature of a dual hyperbolic curve $\tilde{e}_1(\bar{s}_1)$ (timelike ruled surface $\varphi_{e_1}$) as follows:

$$\begin{cases} \bar{s}_1 = \int_0^{s_1} \left\| \tilde{e}_1'(u_1) \right\| du_1 = \int_0^{s_1} (1-\varepsilon \Delta_1) du_1 = s_1 - \varepsilon \int_0^{s_1} \Delta_1 du_1, \quad \Delta_1 = \det(\vec{c}_1', \vec{e}_1, \vec{t}_1), \\ \delta_1 = \langle \vec{c}_1', \vec{e}_1 \rangle, \quad \gamma_1 = -\langle \vec{g}_1', \vec{t}_1 \rangle, \quad \bar{\gamma}_1 = \gamma_1 + \varepsilon(\delta_1 + \gamma_1 \Delta_1) \end{cases} \tag{32}$$

$$\begin{cases} \cosh \bar{\rho}_1 = -\frac{\bar{\gamma}_1}{\sqrt{\left|1-\bar{\gamma}_1^2\right|}}, \quad \sinh \bar{\rho}_1 = -\frac{1}{\sqrt{\left|1-\bar{\gamma}_1^2\right|}}, \quad \text{if } |\bar{\gamma}_1| > 1. \\ \sinh \bar{\rho}_1 = -\frac{\bar{\gamma}_1}{\sqrt{\left|1-\bar{\gamma}_1^2\right|}}, \quad \cosh \bar{\rho}_1 = -\frac{1}{\sqrt{\left|1-\bar{\gamma}_1^2\right|}}, \quad \text{if } |\bar{\gamma}_1| < 1. \end{cases} \tag{33}$$

and

$$\bar{R}_1 = \begin{cases} -\sinh \bar{\rho}_1, & \text{if } |\bar{\gamma}_1| > 1, \\ -\cosh \bar{\rho}_1, & \text{if } |\bar{\gamma}_1| < 1. \end{cases} \quad \text{and} \quad \bar{\gamma} = \begin{cases} \coth \bar{\rho}, & \text{if } |\bar{\gamma}| > 1, \\ \tanh \bar{\rho}, & \text{if } |\bar{\gamma}| < 1. \end{cases} \tag{34}$$

## 5. Darboux Approach to Mannheim Offsets of Spacelike Ruled Surfaces

Let $\varphi_e$ be a spacelike ruled surface generated by dual spacelike unit vector $\tilde{e}$ and $\varphi_{e_1}$ be a ruled surface generated by dual unit vector $\tilde{e}_1$ and let $\{\tilde{e}(\bar{s}), \tilde{t}(\bar{s}), \tilde{g}(\bar{s})\}$ and $\{\tilde{e}_1(\bar{s}_1), \tilde{t}_1(\bar{s}_1), \tilde{g}_1(\bar{s}_1)\}$ denote the dual Darboux frames of $\varphi_e$ and $\varphi_{e_1}$, respectively. Then $\varphi_e$ and $\varphi_{e_1}$ are called Mannheim surface offsets, if

$$\tilde{g}(\bar{s}) = \tilde{t}_1(\bar{s}_1) \tag{35}$$

holds, where $\bar{s}$ and $\bar{s}_1$ are the dual arc-lengths of $\varphi_e$ and $\varphi_{e_1}$, respectively. By this definition, it is seen that the Mannheim offset of $\varphi_e$ is a timelike ruled surface, but the generator of this surface can be timelike or spacelike. In this study, we consider the Mannheim offset surface with timelike ruling. If one considers the second situation, ruling is spacelike, and uses the Definition 3.1 (ii), similar results can be found. Thus, in this paper when we talk about the surfaces $\varphi_e$ and $\varphi_{e_1}$, we mean that $\varphi_e$ is a spacelike ruled surface and $\varphi_{e_1}$ is a timelike ruled surface with timelike ruling and for short we don't write the Lorentzian characters of the surfaces hereinafter.

Let now the ruled surfaces $\varphi_e$ and $\varphi_{e_1}$ form a Mannheim surface offset. Then by considering (35) the relationship between the trihedrons of the surfaces $\varphi_e$ and $\varphi_{e_1}$ can be given as follows

$$\begin{pmatrix} \tilde{e}_1 \\ \tilde{t}_1 \\ \tilde{g}_1 \end{pmatrix} = \begin{pmatrix} \sinh\bar{\theta} & \cosh\bar{\theta} & 0 \\ 0 & 0 & 1 \\ \cosh\bar{\theta} & \sinh\bar{\theta} & 0 \end{pmatrix} \begin{pmatrix} \tilde{e} \\ \tilde{t} \\ \tilde{g} \end{pmatrix} \qquad (36)$$

where $\bar{\theta} = \theta + \varepsilon\theta^*$, $(\theta, \theta^* \in IR)$ is the dual Lorentzian angle between the dual generators $\tilde{e}$ and $\tilde{e}_1$ of Mannheim ruled surfaces $\varphi_e$ and $\varphi_{e_1}$. The real angle $\theta$ is called the offset angle which is the angle between the rulings $\vec{e}$ and $\vec{e}_1$, and $\theta^*$ is called the offset distance which is measured from the striction point $\vec{c}$ of $\varphi_e$ to striction point $\vec{c}_1$ of $\varphi_{e_1}$. And from (36) we may write $\vec{c}_1 = \vec{c} + \theta^*\vec{t}$. Then, $\bar{\theta} = \theta + \varepsilon\theta^*$ is called dual Lorentzian offset angle of the Mannheim ruled surfaces $\varphi_e$ and $\varphi_{e_1}$. If $\theta = 0$, then the Mannheim surface offsets are said to be right offsets.

Now, we give some theorems and results characterizing Mannheim surface offsets.

**Theorem 5.1.** *Let $\varphi_e$ and $\varphi_{e_1}$ form a Mannheim surface offset. The offset angle $\theta$ and the offset distance $\theta^*$ are given by*

$$\theta = -s + c, \quad \theta^* = -\int_0^s \Delta du + c^* \qquad (37)$$

*respectively, where $c$ and $c^*$ are real constants.*

**Proof.** Suppose that $\varphi_e$ and $\varphi_{e_1}$ form a Mannheim offset. Then from (36) we have

$$\tilde{e}_1 = \sinh\bar{\theta}\,\tilde{e} + \cosh\bar{\theta}\,\tilde{t}. \qquad (38)$$

By differentiating (38) with respect to $\bar{s}$ we get

$$\frac{d\tilde{e}_1}{d\bar{s}} = \cosh\bar{\theta}\left(1 + \frac{d\bar{\theta}}{d\bar{s}}\right)\tilde{e} + \sinh\bar{\theta}\left(1 + \frac{d\bar{\theta}}{d\bar{s}}\right)\tilde{t} + \bar{\gamma}\cosh\bar{\theta}\,\tilde{g}. \qquad (39)$$

Since $\varphi_e$ and $\varphi_{e_1}$ form a Mannheim offset, $\dfrac{d\tilde{e}_1}{d\bar{s}}$ and $\tilde{g}$ are linearly dependent. Then, from (39) we get $\dfrac{d\bar{\theta}}{d\bar{s}} = -1$ and for the dual constant $\bar{c} = c + \varepsilon c^*$ we write

$$d\bar{\theta} = -d\bar{s}$$
$$\bar{\theta} = -\bar{s} + \bar{c}$$
$$\theta + \varepsilon\theta^* = -s - \varepsilon s^* + c + \varepsilon c^*$$

Then from (17) we have

$$\theta = -s + c, \quad \theta^* = -\int_0^s \Delta du + c^*,$$

where $c$ and $c^*$ are real constants.

From (37), the following corollary can be given.

**Corollary 5.1.** *Let $\varphi_e$ and $\varphi_{e_1}$ form a Mannheim surface offset. Then $\varphi_e$ is developable if and only if offset distance is constant, i.e. $\theta^* = c^* = $ constant.*

**Theorem 5.2.** *Let $\varphi_e$ and $\varphi_{e_1}$ form a Mannheim surface offset. Then there is the following differential relationship between the dual arc-length parameters of $\varphi_e$ and $\varphi_{e_1}$*

$$\frac{d\bar{s}_1}{d\bar{s}} = \bar{\gamma}\cosh\bar{\theta}. \tag{40}$$

**Proof.** Suppose that $\varphi_e$ and $\varphi_{e_1}$ form a Mannheim offset. Then, from Theorem 5.1 we have

$$\frac{d\tilde{e}_1}{d\bar{s}_1} = \tilde{t}_1 = \bar{\gamma}\cosh\bar{\theta}\frac{d\bar{s}}{d\bar{s}_1}\tilde{g}. \tag{41}$$

Since $\varphi_e$ and $\varphi_{e_1}$ form a Mannheim offset we have $\tilde{t}_1 = \tilde{g}$. Then (41) gives us

$$\bar{\gamma}\cosh\bar{\theta}\frac{d\bar{s}}{d\bar{s}_1} = 1 \tag{42}$$

and from (42) we get (40).

**Corollary 5.2.** *Let $\varphi_e$ and $\varphi_{e_1}$ form a Mannheim surface offset. Then there are the following relationships between the real arc-length parameters of $\varphi_e$ and $\varphi_{e_1}$*

$$\frac{ds_1}{ds} = \gamma\cosh\theta, \quad \frac{dsds_1^* - ds^*ds_1}{ds^2} = \theta^*\gamma\sinh\theta - (\delta + \gamma\Delta)\cosh\theta. \tag{43}$$

**Proof.** Let $\varphi_e$ and $\varphi_{e_1}$ form a Mannheim surface offset. Then by Theorem 5.2, (40) holds. By considering (4), the real and dual parts of (40) are

$$\frac{ds_1}{ds} = \gamma\cosh\theta, \quad \frac{dsds_1^* - ds^*ds_1}{ds^2} = \theta^*\gamma\sinh\theta - (\delta + \gamma\Delta)\cosh\theta \tag{44}$$

which are desired equalities.

In Corollary 5.1, we give the relationship between the offset distance $\theta^*$ and developable spacelike ruled surface $\varphi_e$. Now we give the condition for $\varphi_{e_1}$ to be developable according to $\theta^*$.

From (17) and (32) we have

$$ds^* = \Delta ds, \quad ds_1^* = -\Delta_1 ds_1, \tag{45}$$

respectively. Then writing (45) in (44) and using (43) we get

$$\Delta_1 = -\theta^*\tanh\theta + \frac{\delta}{\gamma}$$

and give the following corollaries:

**Corollary 5.3.** *Let $\varphi_e$ and $\varphi_{e_1}$ form a Mannheim surface offset. Then*

$$\Delta_1 = -\theta^*\tanh\theta + \frac{\delta}{\gamma} \tag{46}$$

*holds.*

**Corollary 5.4.** *Let $\varphi_e$ and $\varphi_{e_1}$ form a Mannheim surface offset. Then $\varphi_{e_1}$ is developable if and only if $\theta^* = \frac{\delta}{\gamma}\coth\theta$ holds.*

**Theorem 5.3.** *Let $\varphi_e$ and $\varphi_{e_1}$ form a Mannheim surface offset. There exists the following relationship between the invariants of the surfaces and offset angle $\theta$, offset distance $\theta^*$,*

$$\delta_1 = \frac{\delta}{\gamma}\tanh\theta - \theta^*. \tag{47}$$

**Proof.** Let the striction lines of $\varphi_e$ and $\varphi_{e_1}$ be $c(s)$ and $c_1(s_1)$, respectively, and let $\varphi_e$ and $\varphi_{e_1}$ form a Mannheim surface offset. Then, from the Mannheim condition we can write

$$\vec{c}_1 = \vec{c} + \theta^*\vec{g}. \tag{48}$$

Differentiating (48) with respect to $s_1$ we have

$$\frac{d\vec{c}_1}{ds_1} = \left(\frac{d\vec{c}}{ds} + \theta^*\gamma\vec{t} + \frac{d\theta^*}{ds}\vec{g}\right)\frac{ds}{ds_1}. \tag{49}$$

From (32) we know that $\delta_1 = \langle d\vec{c}_1/ds_1, e_1\rangle$. Then from (36) and (49) we obtain

$$\delta_1 = \left(\sinh\theta\langle d\vec{c}/ds, \vec{e}\rangle + \cosh\theta\langle d\vec{c}/ds, \vec{t}\rangle + \theta^*\gamma\cosh\langle\vec{t},\vec{t}\rangle\right)\frac{ds}{ds_1}. \tag{50}$$

Since $\langle d\vec{c}/ds, e\rangle = \delta$, $\langle d\vec{c}/ds, \vec{t}\rangle = 0$, $\langle \vec{t}, \vec{t}\rangle = -1$, from (50) we write

$$\delta_1 = \left(\delta\sinh\theta - \theta^*\gamma\cosh\theta\right)\frac{ds}{ds_1}. \tag{51}$$

Furthermore, from (43) we have

$$\frac{ds}{ds_1} = \frac{1}{\gamma\cosh\theta}. \tag{52}$$

Then substituting (52) in (51) we obtain

$$\delta_1 = \frac{\delta}{\gamma}\tanh\theta - \theta^*.$$

**Theorem 5.6.** *If $\varphi_e$ and $\varphi_{e_1}$ form a Mannheim surface offset, then the relationship between the conical curvature $\gamma_1$ of $\varphi_{e_1}$ and offset angle $\theta$ is given by*

$$\gamma_1 = -\tanh\theta. \tag{53}$$

**Proof.** From (32) and (36) we have

$$\begin{aligned}\gamma_1 &= -\langle\vec{g}_1', \vec{t}_1\rangle \\ &= -\left\langle\frac{d}{ds_1}(\cosh\theta\vec{e} + \sinh\theta\vec{t}), \vec{g}\right\rangle \\ &= -\gamma\sinh\theta\frac{ds}{ds_1}\end{aligned} \tag{54}$$

From the first equality of (43) and (54), we have $\gamma_1 = -\tanh\theta$.

**Theorem 5.7.** *If the surfaces $\varphi_e$ and $\varphi_{e_1}$ form a Mannheim surface offset, then the dual conical curvature $\bar{\gamma}_1$ of $\varphi_{e_1}$ is obtained as*

$$\bar{\gamma}_1 = -\tanh\bar{\theta}. \tag{62}$$

**Proof.** From (32), (47), (53) and (61) by direct calculation we have (62).

From (62) we have the following corollary.

**Theorem 5.8.** *If the surfaces $\varphi_e$ and $\varphi_{e_1}$ form a Mannheim surface offset, then the dual curvature of $\varphi_{e_1}$ is given by*

$$\overline{R}_1 = \cosh\overline{\theta}. \tag{63}$$

**Proof.** From (62) we have

$$\sqrt{\left|1-\overline{\gamma}_1^2\right|} = \operatorname{sech}\theta - \varepsilon\theta^* \tanh\theta \operatorname{sech}\theta = \operatorname{sech}\overline{\theta}.$$

Then from (34) we have

$$\overline{R}_1 = \frac{1}{\sqrt{\left|1-\overline{\gamma}_1^2\right|}} = \frac{\operatorname{sech}\theta + \varepsilon\theta^* \tanh\theta \operatorname{sech}\theta}{\operatorname{sech}^2\theta}$$

$$= \cosh\theta + \varepsilon\theta^* \sinh\theta$$
$$= \cosh\overline{\theta}$$

Then we can give the following corollaries.

**Corollary 5.9.** *If $\varphi_e$ and $\varphi_{e_1}$ form a Mannheim surface offset and $|\overline{\gamma}_1|<1$, then the dual spherical radius of curvature of $\varphi_{e_1}$ is given by*

$$\cosh\overline{\rho}_1 = -\cosh\overline{\theta}. \tag{64}$$

**Corollary 5.11.** *Let $\varphi_e$ and $\varphi_{e_1}$ form a Mannheim surface offset. Then $\rho_1^* = \pm\theta^*$ holds.*

If we assume that $|\overline{\gamma}_1|>1$, then we have equalities for a timelike ruled surface whose Darboux vector is timelike and the obtained equalities will be analogue to given ones.

## 6. Conclusions
The dual geodesic trihedron(dual Darboux frame) of a spacelike ruled surface is introduced. Then the characterizations of Mannheim offsets of spacelike ruled surfaces are given in view of dual Darboux frame and new relations between the invariants of Mannheim surface offsets are obtained. The results of the paper are new characterizations of Mannheim surface offsets and also give the relationships for Mannheim surface offsets to be developable according to offset angle and offset distance. Moreover, the relationship between the developable Mannheim surface offsets and their striction lines is given.

**References**
[1] Abdel-All N. H., Abdel-Baky R. A., Hamdoon F. M., *Ruled surfaces with timelike rulings,* App. Math. and Comp., 147 (2004) 241–253.
[2] Blaschke, W., Differential Geometrie and Geometrischke Grundlagen ven Einsteins Relativitasttheorie Dover, New York, (1945).
[3] Dimentberg, F. M., The Screw Calculus and its Applications in Mechanics, (Izdat. Nauka, Moscow, USSR, 1965) English translation: AD680993, Clearinghouse for Federal and Scientific Technical Information.
[4] Hacısalihoğlu. H.H., Hareket Geometrisi ve Kuaterniyonlar Teorisi, Gazi Üniversitesi Fen-Edb. Fakültesi, (1983).
[5] Hoschek, J., Lasser, D., Fundamentals of computer aided geometric design, Wellesley, MA:AK Peters, (1993).

[6] Karger, A., Novak, J., Space Kinematics and Lie Groups, STNL Publishers of Technical Lit., Prague, Czechoslovakia (1978).

[7] Küçük, A., Gürsoy O., *On the invariants of Bertrand trajectory surface offsets,* App. Math. and Comp., 151 (2004) 763-773.

[8] O'Neill, B., Semi-Riemannian Geometry with Applications to Relativity, Academic Press, London, (1983).

[9] Orbay, K., Kasap, E., Aydemir, I, *Mannheim Offsets of Ruled Surfaces,* Mathematical Problems in Engineering, Volume 2009, Article ID 160917.

[10] Önder, M., Uğurlu, H.H., *On the Developable of Mannheim offsets of timelike ruled surfaces in Minkowski 3-space*, arXiv:0906.2077v5 [math.DG].

[11] Önder, M., Uğurlu, H.H., *Some Results and Characterizations for Mannheim Offsets of the Ruled Surfaces*, arXiv:1005.2570v3 [math.DG].

[12] Önder, M., Uğurlu, H.H., *Mannheim Offsets of the Timelike Ruled Surfaces with Spacelike Rulings in Dual Lorentzian Space,* arXiv:1007.2041v2 [math.DG]. arXiv:1005.2570v3 [math.DG].

[13] Önder, M., Uğurlu, H.H., *Dual Darboux Frame of a Timelike Ruled Surface and Darboux Approach to Mannheim Offsets of Timelike Ruled Surfaces,* arXiv:1108.6076v2 [math.DG].

[14] Önder, M., Uğurlu, H.H., Kazaz, M., *Mannheim offsets of spacelike ruled surfaces in Minkowski 3-space,* arXiv:0906.4660v3 [math.DG].

[15] Papaioannou, S.G., Kiritsis, D., *An application of Bertrand curves and surfaces to CAD/CAM,* Computer Aided Design, 17 (8) (1985) 348-352.

[16] Pottmann, H., Lü, W., Ravani, B., *Rational ruled surfaces and their offsets,* Graphical Models and Image Processing, 58 (6) (1996) 544-552.

[17] Ravani, B., Ku, T. S., *Bertrand Offsets of ruled and developable surfaces,* Comp. Aided Geom. Design, 23 (2) (1991) 145-152.

[18] Uğurlu, H.H., Çalışkan, A,. *The Study Mapping for Directed Spacelike and Timelike Lines in Minkowski 3-Space $IR_1^3$,* Mathematical and Computational Applications, 1 (2) (1996) 142-148.

[19] Uğurlu, H. H., Önder, M., *On Frenet Frames and Frenet Invariants of Skew Spacelike Ruled Surfaces in Minkowski 3-space,* VII. Geometry Symposium, 07-10 July 2009, Kırşehir, Turkey.

[20] Uğurlu, H. H., *The relations among instantaneous velocities of trihedrons depending on a spacelike ruled surface,* Hadronic Journal, 22 (1999) 145-155.

[21] Veldkamp, G.R., *On the use of dual numbers, vectors and matrices in instantaneous spatial kinematics,* Mechanism and Machine Theory, II (1976) 141-156.

[22] Wang, F., Liu, H., *Mannheim partner curves in 3-Euclidean space,* Mathematics in Practice and Theory, 37 (1) (2007) 141-143.